\documentclass[fleqn]{article}
\usepackage{amsmath}
\usepackage[rightcaption]{sidecap}
\usepackage{amssymb}
\usepackage{amsthm}
\usepackage{newlfont}
\usepackage[nice]{nicefrac}

\usepackage{ifpdf}
\ifpdf
\usepackage[pdftex]{graphicx}
\else
\usepackage[dvips]{graphicx}
\fi

\ifpdf
\usepackage[pdftex,colorlinks,linkcolor=blue,citecolor=blue]{hyperref}
\else
\usepackage{hyperref}
\fi

\begin{document}


\title{On the expansion of the giant component in percolated $(n,d,\lambda)$ graphs\\
\date{Department of Computer Science and Applied Mathematics\\
 Weizmann Institute, Rehovot 76100, Israel\\
{\tt eran.ofek@weizmann.ac.il} \\
 \vspace{2ex} \today}}

\author{Eran Ofek}
\maketitle

\newtheorem{thm}{Theorem}[section]
\newtheorem{corollary}[thm]{Corollary}
\newtheorem{lemma}[thm]{Lemma}
\newtheorem{prop}[thm]{Proposition}
\newtheorem{claim}[thm]{Claim}

\newtheorem{conj}[thm]{Conjecture}

\newcommand{\bits}{\{ 0,1 \}}
\newcommand{\e}{\epsilon}
\newcommand{\norm}[2]{\| #2 \|_{_{#1}}}
\newcommand{\equivalent}{\stackrel{\rm [1]}{\approx}}
\newcommand{\eqdef}{\stackrel{\vartriangle}{=}}
\newcommand{\Ex}[2]{\underset{#1}{\mathbb E} [#2] }
\newcommand{\discup}{$\bigcup$ \hspace{-2.7ex} $\cdot$ \hspace{0.5ex} }
\newcommand{\twolines}[2]{\stackrel{#1}{#2}}
\newcommand{\pfrac}[2]{\left(\frac{#1}{#2}\right)}
\newcommand{\vs}{} 
\newcommand{\svs}{} 


\def\draft{0}   

\ifnum\draft=1 
    \def\ShowAuthNotes{1}
\else
    \def\ShowAuthNotes{0}
\fi

\ifnum\ShowAuthNotes=1
  \newcommand{\authnote}[2]{{ \bf [#1's note: #2]}}
\else
  \newcommand{\authnote}[2]{}
\fi

\begin{abstract}
Let $d \geq d_0$ be a sufficiently large constant. A $(n,d,c
\sqrt{d})$ graph $G$ is a $d$-regular graph over $n$ vertices whose
second largest (in absolute value) eigenvalue is at most $c
\sqrt{d}$. For any $0 < p < 1, ~G_p$ is the graph induced by
retaining each edge of $G$ with probability $p$. It is known that
for $p > \frac{1}{d}$ the graph $G_p$ almost surely contains a
unique giant component (a connected component with linear number
vertices). We show that for $p \geq \frac{5c}{\sqrt{d}}$ the giant
component of $G_p$ almost surely has an edge expansion of at least
$\frac{1}{\log_2 n}$.
\end{abstract}


\section{Introduction}\label{sec: introduction}
This paper deals with the affect of percolation on the edge
expansion property of \emph{algebraic expander} graphs. These are
$d$-regular graphs in which the second largest eigenvalue (in
absolute value) $\lambda$ of their adjacency matrix is smaller than
$d/5$. We call such a graph a $(n,d,\lambda)$ algebraic expander. A
more intuitive (combinatorial) notion of expansion for a finite
graph $G$ is the \emph{edge expansion}, defined as:
$$
c_E(G) = \inf_{S \subset V_G, |S| < |G|/2} \frac{|\partial_E S|}
{|S|},
$$ where $\partial_E S$ denotes the set of edges with exactly one
vertex in $S$. It is known (due to Tanner, Alon and Milman
\cite{am}, \cite{tanner}) that algebraic expansion implies also a
lower bound on the edge expansion: for a $(n,d,\lambda)$ algebraic
expander it holds that $c_E(G) \geq \frac{d-\lambda}{2}$. (There is
also an inequality in the opposite direction: $c_E(G) \leq
\sqrt{2d(d- \lambda)}$, see \cite{a86a} for details).

Expander graphs received a considerable amount of attention in the
literature in recent years, mostly because these graphs have
numerous applications in theoretical computer science; see, for
example, \cite{AlonSp02, inw, ss, r}. It is well known that for any
fixed $d \geq 3$, random $d$-regular graphs of size $n$ are
asymptotically almost surely expanders, as $n$ grows. The problem of
constructing infinite families of bounded degree expanders is more
difficult, and there are several known constructions of this type
\cite{m, lps, rvw, bl}. The result in this paper applies to the
constructions of \cite{lps,rvw,bl}.

Various applications of expanders rely on their fault-tolerance as
networks. For example, after deleting an appropriate constant
fraction of the edges (arbitrarily), the remaining graph still
contains some linear size connected components or some linear size
paths; see \cite{ac,upfal}. We show that for algebraic expanders if
the deletions are random and independent then with high probability
(with probability that tends to $1$ as $n$ increases) the giant
component has an edge expansion proportional to $\frac{c
\sqrt{d}}{\log_2 n}$. Up to constants, this bound is tight since
with probability bounded away from $0$, the giant component will
contain a $\frac{\log_2 n}{2c \sqrt{d}}$ long "chain" of vertices
each of them, except the first and the last ones, has degree of
exactly $2$ in the giant component. The edge expansion of such a
"chain" is $\frac{4c \sqrt{d}}{\log_2 n}$.

Given a graph $G$, we use $G_p$ to denote the subgraph of $G$
obtained by retaining each edge of $G$ independently with
probability $p$. The graph $G_p$ is the \emph{percolated} version
of $G$. For any graph property of $G$ one can ask if this property
is almost surely retained in $G_p$. A well studied example is the
existence and the uniqueness of a giant component. Roughly
speaking, a giant component is a connected component of $G_p$ that
contains linear fraction of vertices. A question of the same
flavour can be asked also for an infinite graph $G$: for which
values of $p$, $G_p$ is likely to contain an infinite cluster
(connected component) ? is the infinite cluster likely to be
unique ? For several types of graphs, e.g. the $d$ dimensional
grid, the finite/infinite versions turned out to be related. For
many interesting graphs the probability of containing a giant
component (or infinite cluster in the infinite case) exhibits a
sharp threshold around some value called the \emph{critical
probability} (this is due to $0/1$ laws). The critical probability
is denoted by $p_c$. For values of $p$ slightly smaller than $p_c$
the probability for giant component is close to $0$ and for $p$
slightly larger than $p_c$ the probability for giant component is
close to $1$. Benjamini and Schramm \cite{bs} showed that if $G$
is an infinite graph with a positive vertex Cheeger constant
$c_V(G) > 0$ (the Cheeger constant can be defined with respect to
the vertex boundary), then the critical probability for the
existence of an infinite cluster in $G_p$ is $< \frac{1}{1+
c_V(G)} <1$. They also observed that their proof can be applied to
the finite case. Their technique can be easily applied also to the
edge Cheeger constant as shown in \cite{ps}.

A family of expanders is a sequence of $d$-regular graphs $G(n)$,
where $G(n)$ has $n$ vertices and edge expansion of least $b>0$
(independent of $n$). Alon, Benjamini and Stacey \cite{abs} studied
the existence and uniqueness of a giant component when percolation
is applied to families of edge expander graphs. One of their results
is about expander families with increasing girth (the girth of a
graph $G$ is the length of minimum size cycle in it). They show that
for an expander family $G(n)$, with $\text{girth}(G(n))$ that goes
to infinity as $n$ increases, the critical probability $p_c$ for the
existence (and uniqueness) of a giant component is exactly
$\frac{1}{d-1}$. Specifically, for any fixed $\e$, and $p \geq
\frac{1 + \e}{d-1}$ w.h.p. (with high probability, i.e. with
probability that goes to $1$ as $n$, the size of graph, goes to
infinity) $G_p$ contains a connected component with a linear number
of vertices. The fraction of vertices in the giant component depends
on $\e$. The girth, the edge expansion, $d$ and $\e$ influence the
speed in which the probability for a g.c. (giant component) goes to
$1$. For $p \leq \frac{1 - \e}{d-1}$, w.h.p. $G_p$ breaks into
connected components of sub-linear size. It is further shown in
\cite{abs} that if $G(n)$ is an infinite family of $d$-regular
graphs, each one with edge expansion of at least $b >0$, then for
$p$ sufficiently close to $1$ (which depends on $b$) $G(n)_p$ is
w.h.p. a $\frac{1}{\log_2 n}$ expander. They leave as an open
problem the values of $p$ which are close (from above) to the
critical probability $p_c$. Notice that $p_c$ can be as small as
$\frac{1}{d-1}$ as in the case of an infinite family of expanders
with girth that goes to infinity.

Percolation of $(n,d,\lambda)$ graphs has been previously studied by
Frieze, Krivelevich and Martin \cite{fkm}. They gave tight results
about the existence and the uniqueness of the giant component when
$\lambda = o(d)$. Specifically, for $p < \frac{1}{d}$ the graph
$G_p$ almost surely contains only connected components of size
$O(\log n)$. For $p > \frac{1}{d}$ the graph $G_p$ has almost surely
a unique giant component and all other components are of size at
most $O(\log n)$.

%
%

\subsection{Our result}

\begin{thm}\label{prop: main proposition}

Let $d \geq d_0$ be a fixed constant, let $G$ be a $(n,d,c
\sqrt{d})$ algebraic-expander and let $p \geq \frac{5 c}{\sqrt{d}}$
(assuming $c < \frac{\sqrt{d}}{5}$). W.h.p.
the edge expansion of the giant component in $G_p$ is at least
$\frac{c \sqrt{d}}{61 \log n}$.
\end{thm}

Theorem \ref{prop: main proposition} implies that in the case of
algebraic expanders even when $p << 1$ the giant component has edge
expansion $\geq \frac{1}{\log_2 n}$. In contrast, the result in
\cite{abs} is based on a weaker assumption (edge expansion greater
than $\e$) but it implies that the giant component has edge
expansion $\geq \frac{1}{\log_2 n}$ only for values of $p$ close to
$1$. While Theorem \ref{prop: main proposition} requires a somewhat
stronger assumption (spectral gap) from $G$, it implies that the
giant component in $G_p$ has expansion $\geq \frac{1}{\log_2 n}$
also for values of $p$ close to $0$ (depending on the degree $d$ and
$\lambda = c\sqrt{d}$).

The main idea in the proof of Theorem \ref{prop: main proposition}
is to iteratively remove from $G_p$ vertices of low degree until we
are left with an induced subgraph $G_p^k$ that has minimal degree
$\geq \frac{3pd}{5}$. Using known techniques it can be shown that
for large enough $d$ this process removes only small fraction of the
vertices. Moreover, the obtained subgraph $G_p^k$ has edge expansion
bounded away from $0$. To show that the giant component of $G_p$ has
expansion $\geq \frac{1}{\log_2 n}$ (which is best possible up to
constants) we need to handle sets of the giant component that
contain vertices from $OUT \eqdef V \setminus G_p^k$. To do this it
is enough to show that in the graph induced by $G_p$ on $OUT$, the
connected components are smaller than $\log_2 n$. Following the work
of Alon and Kahale \cite{ak} on coloring random $3$-colorable
graphs, several papers \cite{ChenFr96,Flaxman03,GoerdtLa04,Coja05}
dealt with similar versions of this problem: proving that a set
$OUT$ which is the outcome of some procedure applied to a random
graph has no large connected components. Yet, in all the above cases
the graph model was a simple variant of the $G_{n,p}$ model. Our
result can be though of as a "derandomization" of the previous
results as we deal with predetermined constant degree
"pseudo-random" graphs for which there is less randomness in the
induced model.

\medskip

\textbf{Remarks:}
\begin{enumerate} \item Possibly, Theorem
\ref{prop: main proposition} can be extended also for values of $p >
\frac{c}{\sqrt{d}}$, using the same proof technique. To keep the
proof simple, we did not try to optimize this constant.
\item
Theorem \ref{prop: main proposition} holds also for $d$ which is a
function of $n$ if one of the following holds: $c\sqrt{d} = \lambda
= \omega(\log d$) or $d = o(n)$.
\end{enumerate}

%
%

\subsection{Notation}
For a set $U \subset V$, $G[U]$ denotes the subgraph induced by the
edges of $G$ on the vertices of $U$. We use $e(S,S)$ to denote twice
the number of edges having only vertices in $S$. The graph induced
by retaining each edge of $G$ independently with probability of $p$
is denoted by $G_p$. The degree of a vertex $v$ inside a graph $G$
is denoted by $\text{deg}^G_v$. The second largest eigenvalue (in
absolute value) of $G$ is denoted by $\lambda$. We use the term with
high probability (w.h.p) to denote a sequence of probabilities that
tends to $1$ as $n$, the size of $G$, goes to infinity.

\subsection{Spectral gap and pseudo-randomness}
In the following proofs we will use the fact that a graph $G$ with a
noticeable spectral gap is pseudo-random. This is formulated by the
following Lemma also known as the expander mixing lemma (see
\cite{AlonSp02} for proof).

\begin{lemma}
Let $G$ be a $d$-regular graph with second largest (in absolute
value) eigenvalue $\lambda$. Then, for any $S,T \subseteq V$:
\begin{align*}
| e(S,T) -\frac{d}{n} |S||T| ~| < \lambda \sqrt{|S||T|},
\end{align*}
 where $e(S,T)$ is the
number of directed edges from $S$ to $T$ in the adjacency matrix of
$G$.

\end{lemma}

In terms of undirected edges (when $G$ is undirected), $e(S,T)$
equals the number of edges between $S \setminus T$ to $T$ plus twice
the number of edges that contain only vertices of $S \cap T$.

\begin{corollary}\label{cor: small sets are sparse}
Let $G$ be a $(n,d,c\sqrt{d})$ algebraic expander. For any set $U$
of size $\leq \frac{cn}{k \sqrt{d}}$ the average degree in $G[U]$ is
at most $c\sqrt{d}(1 + 1/k)$.
\end{corollary}

\begin{proof}
The number of edges inside $G[U]$ is $e(U,U)/2$ since every edge
whose both endpoints are in $U$ is in fact two directed edges from
$U$ to $U$. It follows that the average degree in $G[U]$ is
$\frac{e(U,U)}{|U|}$. By the expander mixing lemma:
$$ e(U,U) \leq \frac{d|U|^2 }{n}  + c \sqrt{d} |U| \leq  c\sqrt{d}|U|
\left( 1+ \frac{\sqrt{d}|U|}{c n} \right) \leq c\sqrt{d}|U|(1 +
1/k)$$
\end{proof}

We will frequently use the fact that small enough sets in $G$ are
rather sparse, as stated in Corollary \ref{cor: small sets are
sparse}. When $c$ close to its smallest possible value for constant
degree graphs (i.e. $\lambda = c\sqrt{d} \approx 2\sqrt{d-1}$, see
\cite{a86a} for details) there is a slightly stronger bound on the
density of small sets given by \cite{k}. We do not use this stronger
bound as it gives asymptotically the same result for large values of
$c$.


\section{Proof of Theorem \ref{prop: main proposition}}

We use a process similar to \cite{upfal} which is aimed to reveal
a large edge expanding subgraph of $G_p$. Let $S_0 = \{v \in G_p :
\deg_v^{G_p} \not\in [\frac{4pd}{5}, \frac{6pd}{5}] \}$. We begin
by removing from $G_p$ the vertices of $S_0$. The new induced
graph, which we denote by $G_p^0$, may contain vertices whose
degree in $G_p^0$ is $< \frac{4}{5} pd$, because removing $S_0$
affect the degrees of the remaining vertices. We then obtain a
sequence of subgraphs $G_p^i$ by iteratively removing from
$G_p^{i-1}$ any vertex whose degree in $G_p^{i-1}$ is $<
\frac{3}{5} pd$. The process stops once the minimal degree in the
remaining graph $G_p^k$ is at least $\frac{3}{5} pd$. We denote
the set $V \setminus V(G_p^k)$ by $OUT$.\\

\textbf{Remark:} Since $d$ is fixed and $n \rightarrow \infty$ it
holds that the constant $c$ (from $\lambda = c\sqrt{d}$) is at least
$1$ (in fact $c \geq 1-o(1)$ even if $d$ grows with $n$ but it is
$o(n)$). We will use this fact occasionally.

\subsection{Proof overview}
The main idea in the proof of Theorem \ref{prop: main proposition}
is as follows. We remove from $G_p$ low degree vertices until the
induced graph $G_p^k$ has a large enough minimal degree. We first
show that $G_p^k$ itself has edge expansion of at least
$\frac{pd}{13}$ and contains almost all the vertices of $G$ (thus it
must be contained in a giant component); this part of the proof uses
standard techniques. We then show that in $G_p[OUT]$ the largest
connected component is of size at most $\log_2 n$ (this implies the
uniqueness of the giant component). The above two facts imply that
any set that belongs to the giant component and is entirely in
$V(G_p^k)$ or entirely in $OUT$, has expansion $\geq \frac{1}{\log_2
n}$. Using also the property that any vertex of $V(G_p^k)$ has at
most $\frac{3pd}{5}$ neighbors in $OUT$ we then prove that any set
of the giant component with size $\leq n/2$ has edge expansion $\geq
\frac{1}{\log_2 n}$.

\begin{figure}
\centering
\includegraphics[width=13cm,height=7cm]{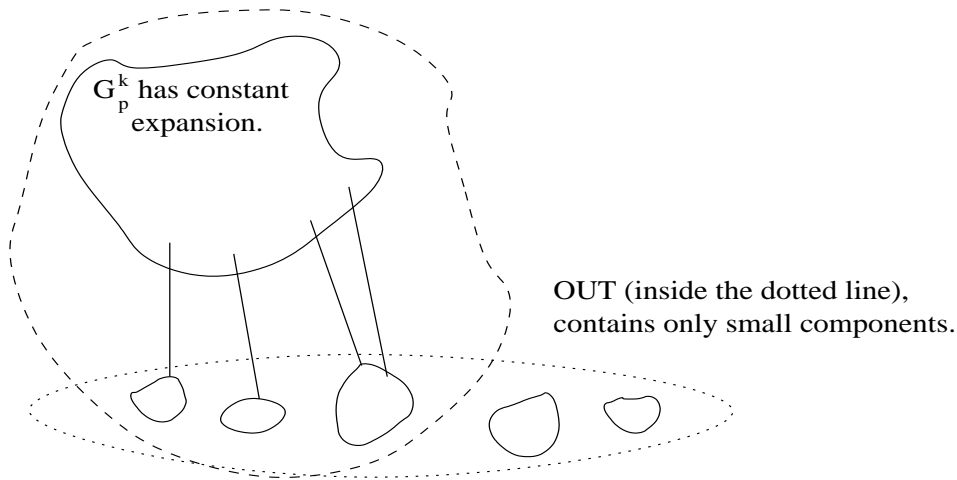}
\caption{The structure of $G_p$.}
\end{figure}

The expected degree in $G_p$, which is $pd$, is large enough so that
only few vertices are removed in the process of extracting $G_p^k$
(namely $OUT$ is small). To show it we use the following idea used
in \cite{ak,upfal}. Initially, the set $S_0$ is small (roughly
$e^{-\Omega(c\sqrt{d})}n$). Every vertex which is removed in the
iteration process has at least $\frac{pd}{5}$ edges to vertices that
were previously removed. Thus if the iteration process is too long,
the set $V \setminus V(G_p^i)$  becomes too dense with contradiction
to Corollary \ref{cor: small sets are sparse}.

\begin{lemma}\label{lemma: out is small}
W.h.p. the number of vertices in $G_p^k$ is at least $(1 -
e^{-\frac{1}{12}c \sqrt{d}})n$.
\end{lemma}
\noindent The proof of Lemma \ref{lemma: out is small} is deferred
to Section \ref{proof: out is small}.

The minimal degree in $G_p^k$ is $\frac{3pd}{5} \geq 3c \sqrt{d}$. A
set $S$ smaller than $\frac{cn}{\sqrt{d}}$ contains at most
$c|S|\sqrt{d}$ internal edges (by Corollary \ref{cor: small sets are
sparse}), thus for such set the expansion is at least $\frac{3pd}{5}
- 2c\sqrt{d} \geq \frac{pd}{5}$. To establish the edge expansion of
larger sets, a standard argument using the Chernoff and union bounds
suffices. In Section \ref{proof: expansion of G_k} we give the full
proof of the following Lemma.

\begin{lemma} \label{lemma: expansion of G_k}
W.h.p. the graph $G_p^k$ has an edge expansion of at least
$\frac{pd}{13}$.
\end{lemma}

We next show that the connected components in $G_p[OUT]$ are of size
at most $\log_2 n$. A direct "brute force" approach using the union
bound over all possible trees of size $\log_2 n$ does not seem to
work here because we don't have a good enough upper bound on the
probability that a fixed tree $T$ is in $OUT$. Notice that we can
not simply claim that every vertex in $OUT$ has a low degree in
$G_p$ (if this were true then probably a simple argument would have
sufficed). It may be the case that a vertex in $OUT$ has high degree
in $G_p$ but it is connected (directly or via other vertices) to
vertices of low degree in $G_p$.

Following the work of Alon and Kahale \cite{ak} on coloring random
$3$-colorable graphs, several papers
\cite{ChenFr96,Flaxman03,GoerdtLa04,Coja05} dealt with similar
versions of this problem (proving that a set $OUT$ which is the
outcome of some procedure applied to a random graph has no large
connected components). Yet, the analysis they gave is rather
complicated. The reason for the difficulty is that $OUT$ is not a
random set independent of $G_p$; its vertices are in fact correlated
and depend on the edges of $G_p$. A new aspect of the current paper
is simplifying the proof that $OUT$ has only small connected
components. This is done using a reduction which will be described
shortly. The outcome of the reduction is that instead of having to
prove that w.h.p. $G_p[OUT]$ has no large trees, we need only to
prove that w.h.p. $G_p[OUT]$ has no large \emph{balanced} trees. A
balanced tree is a tree in which at least $\frac{1}{3}$ of its
vertices are low degree vertices in $G_p$. Proving that $G_p[OUT]$
contains no balanced trees can be done directly by applying the
union bound over all possible sets that may form large trees in
$G_p$.

We will now explain the reduction. We show that w.h.p. any maximal
connected component $U$ in $G_p[OUT]$ is
\emph{$\frac{1}{3}$-balanced} i.e. at least $1/3$ of its vertices
are from $S_0$. The argument is as follows: every vertex removed in
iterations $1,2,..,k$ has at least $\frac{pd}{5}$ edges to
previously removed vertices. Thus, if at least $2/3$ of the vertices
of $U$ are added during the iterations then the average degree in
$G[U]$ is at least $\frac{2}{3} 2 c\sqrt{d}$. On the other hand,
since $|U| < e^{-\frac{1}{12}c\sqrt{d}}n$ by Corollary \ref{cor:
small sets are sparse} the average degree in $G[U]$ is at most
$c\sqrt{d}(1+ e^{-c\sqrt{d}/2}$), which yields a contradiction.
Having established that w.h.p. every maximal connected component in
$G[OUT]$ is balanced, our problem (showing that w.h.p. there are no
large connected components in $G[OUT]$) is reduced to showing that
w.h.p. any maximal balanced connected component of $G[OUT]$ is of
size $< \log_2 n$. Since proving the claim on balanced  trees is
simpler (it is a private case) we translate the claim into a claim
on trees. This is done by showing that any $\frac{1}{3}$-balanced
connected component of $G[U]$ of size $\geq \log_2 n$ contains a
$\frac{1}{3}$-balanced tree whose size is in $[\log_2 n, 2 \log_2
n]$; this follows from following Lemma whose proof is deferred to
Section \ref{sec: lemmas proof}.

\begin{lemma}\label{lemma: tree sample}
Let $G$ be a connected graph whose vertices are partitioned into two
sets: $S$ and $I$. Let $\frac{1}{k}$ be a lower bound on the
fraction of $S$ vertices, where $k$ is an integer. For any $1 \leq t
\leq |V(G)|/2$ there exists a tree whose size is in $[t,2t-1]$ and
at least $\frac{1}{k}$ fraction of its vertices are from $S$.
\end{lemma}

To summarize, in order to show that $G_p[OUT]$ has no connected
components of size $\log_2 n$ we need only to prove the following
Claim.

\begin{claim} \label{claim: out has no large balanced trees}
W.h.p. $G_p[OUT]$ has no balanced trees of size in $[\frac{1}{2}
\log_2 n, \log_2 n]$.
\end{claim}

\begin{proof}
We want to bound the probability that $G_p[OUT]$ contains any
$\frac{1}{3}$-balanced tree of size in $[t/2,t]$ (we fix the
parameter $t$ later). Such a tree is called "bad".

$$ \sum_{\text{$T$ is a tree}} \Pr[T \subset OUT ~\wedge~ \text{T is $\tfrac{1}{3}$-balanced}] $$
The number of trees of size $t$ in a $d$-regular graph $G$ is at
most $n d^{2t}$, since each tree can be uniquely mapped into a
closed path of length $2t$. For each tree of size $t$ there are at
most $2^t$ ways of choosing a subset of size $\geq t/3$. Any fixed
set of size $\geq t/3$ is in $S_0$ with probability of at most
$e^{-\frac{1}{2} \frac{1}{25} pd t/3}$. Thus the probability that
there is a balanced bad tree in $OUT$ is at most:
\begin{multline}
 t n d^{2t} 2^{t} e^{-\frac{1}{300} pd t} \leq \exp{(\log t +\log n
 +2t \log d  + t -\tfrac{1}{300} pd t}) \stackrel{\small \underbrace{pd \geq 5c\sqrt{d}}} {\leq} \\
 \exp{( \log n + 3t \log d - \tfrac{c\sqrt{d} t}{60})}= o(1),
 \end{multline}
 for $t \geq \frac{61 \log n}{ c \sqrt{d}}$ (for fixed $d$ and large enough $n$ it holds that $c>1$, see \cite{a86a}).

\emph{Remark:} Notice that in the last inequality we used $\frac{c
\sqrt {d}}{60} > 3 \log d$. This holds also when $d$ is a function
of $n$, if $n$ is large enough and $d = o(n)$ (it is known that
$c\sqrt{d} = \lambda  \geq \sqrt{d - \frac{d^2}{n-1}}$).
\end{proof}

 Since all
connected components of $G[OUT]$ are of size at most $\frac{61
\log_2 n}{c\sqrt{d}}$, sets from $OUT$ that belong to the giant
component have expansion of at least $\frac{c \sqrt{d}}{61 \log_2
n}$. It remains to handle sets of the giant component that
intersects both $V(G_p^k)$ and $OUT$.
\begin{lemma}
W.h.p. any set $S$ that belongs to the giant component and whose
size is at most $n/2$ has edge expansion of at least
$\frac{c\sqrt{d}}{61\log_2 n}$.
\end{lemma}

\begin{proof}
We already handled sets which are completely in $OUT$ or completely
in $V(G_p^k)$. Let $S$ be a set of the g.c. (giant component) that
intersects both $OUT$ and $V(G_p^k)$. Denote by $\bar{S}$ the
complement of $S$ in the giant component. Denote by $S_1,S_2$ the
intersection of $S$ with $OUT,V(G_p^k)$ respectively. We further
partition $S_1$ into $S_{11},S_{12}$ as follows: $S_{11}$ contains
all the connected components of $G_p[S_1]$ that have at least one
edge into $\bar{S}$ and $S_{12}$ contains all the connected
components of $G_p[OUT]$ that have only edges to $S_2$. It is enough
to show that:
$$|E(S_{11},\bar{S})| \geq \frac{|S_{11}| c\sqrt{d}}{61 \log_2
n}, ~~~~|E(S_{12} \cup S_2, \bar{S})| \geq \frac{|S_{12} \cup
S_2|}{18}~.$$ The first inequality follows immediately from the
definition of $S_{11}$. The second inequality is derived as follows:
$|E(S_{12} \cup S_2, \bar{S})| \geq |E(S_2,\bar{S})| \geq |S_2|
pd/13$. Thus
$$ \frac{|E(S_{12} \cup S_2, \bar{S})|}{|S_{12}| + |S_2|} \geq
\frac{|S_2| pd/13 }{|S_{12}| + |S_2|} \geq
\frac{1}{\frac{13}{pd}(|S_{12}|/|S_2| + 1)} \geq \frac{1}{18},$$
where the last inequality holds because every vertex of $G_p^k$ has
at most $\frac{6pd}{5}$ neighbors in $OUT$.
\end{proof}

\subsection{Proofs of lemmas \ref{lemma: out is small}, \ref{proof: expansion of G_k}, \ref{lemma: tree sample}}\label{sec: lemmas proof}

The proofs of Lemmas \ref{lemma: expansion of G_k} and \ref{lemma:
out is small} are rather standard and are based on the fact that
every small enough set $S$ ($ << \frac{cn}{\sqrt{d}}$) contains at
most $|S| c\sqrt{d}(1+o(1))/2$ internal edges.


\begin{proof}[Proof of Lemma \ref{lemma: out is small}]\label{proof: out is small}

A fixed vertex $v$ belongs to $S_0$ with probability $< e^{-
\frac{1}{2} (\frac{1}{5})^2 pd} \leq e^{-\frac{c\sqrt{d}}{10}}$.
Thus, the expected size of $S_0$ is $e^{-\frac{1}{10} c \sqrt{d}}n$.
With probability of $1 - o(1)$ the cardinality of $S_0$ is at most
$e^{-\frac{1}{12} c \sqrt{d}}n$. We briefly sketch the proof. We use
the edge exposure martingale to prove that $S_0$ is concentrated
around its expectation. We fix some order on the $m=nd/2$ edges of
$G$. Let $X_0,X_1,...,X_{m}$ be the martingale sequence, where $X_i$
is the expectation of $S_0$ after exposing the first $i$ edges of
$G_p$. Notice that $X_0 = \mathbb E_{G_p} [S_0]$. The value of
$X_{m}$ is the value of the random variable $S_0$ where the
probability measure is induced by $G_p$. To use the Azuma inequality
we need to upper bound the martingale difference $|X_{i+1} - X_{i}|$
(for $i=0,..,m-1$). It is known that if $S_0$ satisfies the edge
Lipschitz condition with a constant $\Delta$, then also the
martingale difference is bounded by $\Delta$ (see \cite{AlonSp02}).
It is clear that for a fixed graph $G'$, adding/removing a single
edge can change the value of $S_0$ by at most $2$. By Azuma's
inequality:
$$\Pr[X_{m}
> X_0 +  \lambda ] \leq e^{-\lambda^2/(2m \Delta^2)} .$$
Substituting $\lambda = e^{-\frac{1}{10} c \sqrt{d}}n,\Delta =2,
m=nd/2$ we derive that w.h.p. $|S_0|$ is at most $e^{-\frac{1}{12} c
\sqrt{d}}n$.

We next show that the number of vertices removed after removing
$S_0$ (that is $k$) is at most $|2S_0|$. Every vertex that is
removed in the iterative process has at least $\frac{pd}{5} \geq c
\sqrt{d}$ edges which goto previously removed vertices, because its
degree drops from at least $\frac{4pd}{5}$ (as it does not belong to
$S_0$) down to at most $\frac{3pd}{5}$ (at the point it was
removed). By contradiction, assume that $k \geq 2|S_0|$. Consider
the situation immediately after iteration $i = 2|S_0|$. Denote by
$U$ the set of vertices not in $G_p^i$. The average degree in
$G_p[U]$ is at least $\frac{2}{3} 2c\sqrt{d}$. At this point $|U|
\leq 3 e^{-\frac{1}{12} c \sqrt{d}}n$. We derive a contradiction as
by Corollary \ref{cor: small sets are sparse} the average degree in
$G[U]$ is at most $c \sqrt{d}(1 + e^{-\frac{c \sqrt{d}}{15}})$.
\end{proof}

\begin{proof}[Proof of Lemma \ref{lemma: expansion of G_k}]
\label{proof: expansion of G_k} The proof is divided into two parts.
First consider sets of cardinality $\leq \frac{cn}{ \sqrt{d}}$. Fix
a set $S \subset V$. The edge expansion of $S$ (in $G_p^k$) is at
least:
$$\sum_{v \in S} \deg(v) -
 e(S,S),$$ (remember that $e(S,S)$ is twice the number of edges inside
 $S$). Every vertex $v$ of $G_p^k$ has degree of at least
$\frac{3pd}{5} \geq 3c\sqrt{d}$ in $G_p^k$. By Corollary \ref{cor:
small sets are sparse} $e(S,S)$ in $G$ is at most $|S|c\sqrt{d}(1 +
1) = 2|S|c\sqrt{d}$. It follows that the edge expansion of $S$ is at
least $\frac{3pd}{5} - 2c \sqrt{d} \geq \frac{pd}{5}$.

Consider now a set $S \subset V$ of size $\alpha n$ such that
$\frac{c}{\sqrt{d}} < \alpha \leq \frac{1}{2}$. By the expander
mixing lemma, the number of edges between $S$ and $V \setminus S$ in
$G$ is at least:
\begin{align*}
\alpha (1- \alpha)dn - c \sqrt{d} \sqrt{\alpha(1- \alpha)} n &=
\alpha (1- \alpha)dn \left(1 - \tfrac{c}{\sqrt{d \alpha (1 -
\alpha)}}\right) \\&\geq \alpha (1- \alpha)dn/3.
\end{align*}
The last inequality follows from $\frac{c}{\sqrt{d}} \leq \alpha
\leq \frac{1}{2}$ and $c \leq \frac{\sqrt{d}}{5}$. The number of
edges between $S$ and $V \setminus S$ in $G_p$ is at least $p
\alpha(1-\alpha)dn/6$ with probability of $e^{-\frac{1}{8} p \alpha
(1- \alpha)dn/3}$ (follows from the Chernoff bound). The number of
subsets of size $\alpha n$ is at most $\left( \frac{ne}{\alpha
n}\right)^{\alpha n} \leq e^{\alpha n (1+ \log \frac{1}{\alpha})}$.
Thus the probability for a "bad" set of size $\alpha n$ is at most:
\begin{align*}
\exp \left({\alpha n (1+ \log \tfrac{1}{\alpha}) -\tfrac{1}{8}  p
\alpha (1- \alpha)dn/3}\right) &\leq \exp \left({\alpha n ( 1+ \log
\tfrac{1}{\alpha} ~- \tfrac{5c\sqrt{d}(1 - \alpha)}{ 24})} \right) \\
& \leq  \exp \left({- \alpha c n \sqrt{d}/10} \right).
\end{align*}
The last inequality holds for large enough $d$ (and
$\frac{1}{\alpha} < \sqrt{d}$).
 Summing over all values of $\alpha n$
gives that w.h.p. there is no bad set. In other words, every set of
size in $[\frac{cn}{\sqrt{d}},\frac{n}{2}]$ has an edges expansion
of at least $\frac{pd}{12}$ in $G_p$. Since the number of edges that
contain at least one vertex from $OUT$ is bounded by $ d
e^{-\frac{c\sqrt{d}}{12}}n$, we conclude that any subset $U$ of
$V(G_p^k)$ with size $\geq \frac{cn}{\sqrt{d}}$ has edge expansion
of at least
$$\left(\frac{pd}{12} |U| - d e^{-c \sqrt{d}/12}n \right)/|U|
\geq \frac{pd}{13}.$$

\end{proof}

\begin{proof}[Proof of lemma \ref{lemma: tree sample}]
We use the following well know fact: any tree $T$ contains a
\textit{center} vertex $v$ such that each subtree hanged on $v$
contains strictly less than half of the vertices of $T$.

Let $T$ be an arbitrary spanning tree of $G$, with center $v$. We
proceed by induction on the size of $T$. Consider the subtrees
$T_1,...,T_k$ hanged on $v$. If there exists a subtree $T_j$ with at
least $t$ vertices then also $T \setminus T_j$ has at least $t$
vertices. In at least one of $T_j , T \setminus T_j$ the fraction of
$S$ vertices is at least $\frac{1}{k}$ and the lemma follows by
induction on it. Consider now the case in which all the trees have
less than $t$ vertices. If in some subtree $T_j$ the fraction of $S$
vertices is at most $\frac{1}{k}$, then we remove it and apply
induction to $T \setminus T_j$. The remaining case is that in all
the subtrees the fraction of $S$ vertices is strictly more than
$\frac{1}{k}$. In this case we start adding subtrees to the root $v$
until for the first time the number of vertices is at least $t$. At
this point we have a tree with at most $2t-1$ vertices and the
fraction of $S$ vertices is at least $\frac{1}{k}$. To see that the
fraction of $S$ vertices is at least $\frac{1}{k}$, we only need to
prove that the tree formed by $v$ and the first subtree has
$\frac{1}{k}$ fraction of $S$ vertices. Let $r$ be the number of $S$
vertices in the first subtree and let $b$ be the number of vertices
in it. Since $k$ is integer we have: $\frac{r}{b}
> \frac{1}{k} \Longrightarrow \frac{r}{b+1} \geq \frac{1}{k}$.
\end{proof}

\section{Open problems}
We were able to show that a percolation applied to a family of
$d$-regular expander graphs with eigenvalue gap retains some
expansion properties of the original graphs, even when $p$ is close
to $0$. There are still many open problems, we list here two of
them:

\begin{enumerate}
\item Find other classes of expander families that retain expansion
properties after percolated with values of $p$ close to $0$. For
example, a family of expanders with girth that goes to infinity (for
such a family some result is given at \cite{abs} for $p$ close to
$1$).

\item Is Theorem \ref{prop: main proposition} is tight ? If we
drop the requirement that $d$ is a constant and allow it to be a
function of $n$, then the current proof of Claim \ref{claim: out has
no large balanced trees} breaks done (when $d$ is proportional to
$n$). However it is plausible that a different counting argument may
work; one example where a modified argument works is $K_n$ (we have
a proof for this case). If this is the case then Theorem \ref{prop:
main proposition} is tight (up to constant factors) because for the
complete graph $K_n$ it holds that $d= n-1, c= \frac{1}{\sqrt{n-1}}$
and for $p << \frac{c}{\sqrt{d}}$ the percolated graph is not likely
to contain a giant component. Anyway, for constant $d$ the question
is interesting: there is a gap between the critical probability
$\frac{1}{d}$ for which there is a giant component and the
probability $\frac{5c}{\sqrt{d}}$ for which there is
$\frac{1}{\log_2 n}$ edge expansion.
\end{enumerate}

\section*{Acknowledgments}
I thank Itai Benjamini for suggesting the problem and for useful
discussions.
%
%
%
%
%
%



\begin{thebibliography}{99}

\bibitem{a86a}
N. Alon.
\newblock Eigenvalues and expanders.
\newblock Combinatorica, 6(2):83--96, 1986.


\bibitem{abs}
N. Alon, I. Benjamini and A. Stacey.
\newblock Percolation on finite graphs and isoperimetric inequalities.
\newblock Ann. Probab. 32(3):1727--1745, 2004.

\bibitem{ac}
N. Alon and F. R. K. Chung.
\newblock Explicit construction of linear sized tolerant networks.
\newblock Discrete Math., 72(1-3):15--19, 1988.

\bibitem{ak}
N.~Alon and N.~Kahale.
\newblock A spectral technique for coloring random $3$-colorable graphs.
\newblock {\em SIAM Journal on Computing}, 26(6):1733--1748, 1997.

\bibitem{am}
N.~Alon and V.D. Milman.
\newblock $\lambda_1$, isoperimetric inequalities for graphs and
superconcentrators.
\newblock J. of Combinatorial Theory, Ser. B, 38:73--88,1985.

\bibitem{AlonSp02}
N.~Alon and J.~Spencer.
\newblock {\em The Probabilistic Method}.
\newblock John Wiley and Sons, 2002.

\bibitem{bs}
I. Benjamini and O. Schramm.
\newblock  Percolation beyond $\mathbb{Z}^d$, many questions and a few answers.
\newblock  Electr. Comm. Probab, 1:71–-82, 1996.

\bibitem{bl}
Y. Bilu and N. Linial.
\newblock Lifts, discrepancy and nearly optimal spectral gaps.
\newblock FOCS 2004, 404--412.

\bibitem{ChenFr96}
H.~Chen and A.~Frieze.
\newblock Coloring bipartite hypergraphs.
\newblock IPCO 1996, 345--358.

\bibitem{Coja05}
A.~Coja-Oghlan.
\newblock A spectral heuristic for bisecting random graphs.
\newblock SODA 2005, 850--859.

\bibitem{Flaxman03}
A.~Flaxman.
\newblock A spectral technique for random satisfiable 3cnf formulas.
\newblock SODA 2003, 357--363.

\bibitem{fkm}
A. Frieze, M. Krivelevich and R. Martin.
\newblock Emergence of a giant component in random subgraphs of pseudo-random
graphs.
\newblock Random Structures and Algorithms, 24(1):42--50, 2004.

\bibitem{GoerdtLa04}
A.~Goerdt and A.~Lanka.
\newblock On the hardness and easiness of random 4-sat formulas.
\newblock ISAAC 2004, 470--483.

\bibitem{inw}
R. Impagliazzo, N. Nisan and A. Wigderson.
\newblock Pseudorandomness for network algorithms.
\newblock STOC 1994, 356-364.

\bibitem{k}
N. Kahale.
\newblock Eigenvalues and expansion of regular grpahs.
\newblock J. of the ACM, 42(5):1091--1106,1995.

\bibitem{lps}
A. Lubotzky, R. Phillips and P. Sarnak.
\newblock Ramanujan Graphs, Combinatorica, 8:261--277, 1988.

\bibitem{m}
G.A. Margulis.
\newblock Explicit constructions of concentrators.
\newblock Problemy Pereda\v{c}i Informacii, 9(4):71--80,1973; English transl. in Problems
of Information Transmission, 325 - 332, 1975.

\bibitem{ps}
I. Pak and T. Smirnova-Nagnibeda.
\newblock Uniqueness of percolation on nonamenable Cayley graphs.
\newblock Comptes Rendus Acad. Sci. Paris, Ser. I Math, vol. 330(6): 495--500, 2000.

\bibitem{r}
O. Reingold.
\newblock Undirected ST-Connectivity in Logspace.
\newblock STOC 2005, 376--385.

\bibitem{rvw}
O. Reingold, S. Vadhan and A. Wigderson.
\newblock Entropy Waves, the Zig-Zag graph product and new constant-degree expanders.
\newblock Annals of Mathematics, 155:157--187, 2002.

\bibitem{ss}
M. Sipser and D.A. Spielman.
\newblock Expander codes.
\newblock IEEE Transactions on Information Theory 42(6): 1710-1722,
1996.

\bibitem{tanner}
R.M. Tanner.
\newblock Explicit construction of concentrators from generalized n-gons.
\newblock J. Algebraic Discrete Methods, 5:287--294, 1984.

\bibitem{upfal}
E. Upfal.
\newblock Tolerating linear number of faults in networks of bounded degree.
\newblock PODC 1992, pages 83--89.



\end{thebibliography}
\end{document}